\documentclass[11pt]{article}
\usepackage{amsmath,amssymb,mathrsfs}
\usepackage{amsfonts}
\input xy
\xyoption{all}

\newcommand{\pf}{\noindent{\em Proof.}\ }
\newcommand{\qed}{\hfill $\Box$\\}

\renewcommand{\ni}{\noindent}

\newcommand{\ZZ}{\mathbb{Z}}
\newcommand{\QQ}{\mathbb{Q}}
\newcommand{\RR}{\mathbb{R}}
\newcommand{\CC}{\mathbb{C}}

\newcommand{\GG}{\mathbb{G}}
\renewcommand{\SS}{\mathbb{S}}

\newcommand{\UU}{\mathbb{U}}
\newcommand{\EE}{\mathbb{E}}

\newcommand{\tensor}{\otimes}

\newcommand{\rc}{\subset}

\newcommand{\isom}{\simeq}
\newcommand{\isoto}{\tilde{\to}}

\renewcommand{\projlim}{\underleftarrow{\lim}}

\newcommand{\absfr}[1]{\sigma_{#1}}

\newcommand{\Spec}{\operatorname{Spec}}
\newcommand{\Spf}{\operatorname{Spf}}
\newcommand{\Gal}{\operatorname{Gal}}
\newcommand{\End}{\operatorname{End}}
\newcommand{\Hom}{\operatorname{Hom}}
\newcommand{\rk}{\operatorname{rk}}

\begin{document}

\title{Special lifts of ordinary K3 surfaces and applications
	\footnote{This work is supported in part by Professor N. Yui's
		Discovery Grant from NSERC, Canada.}}
\author{\sc Jeng-Daw Yu}
\date{November 2006}
\maketitle

\begin{abstract}
We study the cohomological properties
of quasi-canonical lifts of an ordinary K3 surface over a finite field.
As applications,
we prove a Torelli type theorem for ordinary K3 surfaces
over finite fields
and establish the Hodge conjecture for arbitrary self products
of certain K3 surfaces over the complex numbers.
\end{abstract}

\begin{center}
{\large 0. Introduction}
\end{center}

Let $X$ be an ordinary K3 surface over a finite field $k$.
Then the deformation functor of $X$ admits a natural group structure
(\cite{Nyg-Tate}, Theorem 1.6).
Under this group structure,
the deformation of $X$ over an artinian local ring
corresponding to the identity element
is called the canonical lift of $X$.
Deformations that corresponds to torsion elements
are called quasi-canonical lifts of $X$.
In \cite{Nyg-Tate},
Nygaard uses these lifts to prove the Tate conjecture
for ordinary K3 surfaces over finite fields.\\

In this paper,
we investigate the cohomological properties
of these special lifts of an ordinary K3 surface $X$
to a complete discrete valuation ring of characteristic zero.
We show that these lifts are characterized by
the existence of a frobenous acting on the Betti cohomology of the lifts.
This property is similar to the well-known fact
in the case of ordinary abelian varieties over a finite field
(\cite{Messing}, Appendix).
We also prove a Torelli type theorem
for ordinary K3 surfaces over a finite field
(cf. \cite{Nyg-Torelli}, Theorem 2.1).
We compute the Hodge group of the Hodge structure
associated to a quasi-canonical lift of $X$.
We then apply the result to establish the Hodge conjecture for
arbitrary self products of certain K3 surfaces.\\

This paper consists of two parts and each part contains two sections.
In Part I, \S 1,
we investigate the basic properties of
the quasi-canonical lifts of an ordinary K3 surface
over a finite field,
and in \S 2 we study the Hodge structures of these special lifts.
In Part II,
we give applications.
In \S 3, we prove a Torelli type theorem for ordinary K3 surfaces.
In \S 4, we show that certain K3 surfaces over characteristic zero
are quasi-canoncal lifts of their reductions.
Using this fact,
we translate the validity of the Tate conjecture
for arbitrary self products of the ordinary reductions
(\cite{Zarhin-Tate}) to the validity of the Hodge conjecture
via the explicit description of the Hodge groups
of quasi-canonical lifts.\\

This paper contains the results of my thesis.
I would like to thank my advisor,
Professor Shing-Tung Yau, for his coherent support and encouragement.
I would like to thank Professor Ching-Li Chai
for discussion of this work and answering my questions.
Thanks are also due to Professor Richard Taylor
for his interests in this work.
The study of Hodge structures is also inspired by
the discussions with Professor Noriko Yui
and I would like to thank her heartily. \\

\begin{center}
{\bf Notations}\\
\end{center}

\ni
$k$, a perfect field of characteristic $p > 0$
with countably many elements\\
$W = W(k)$, the Witt vectors of $k$\\
$K = W \tensor_{\ZZ_p} \QQ_p$, the field of fractions of $W$\\
$\absfr{} : k \to k$, $W \to W$, $K \to K$, the absolute frobenius\\
$\bar{K}$, an algebraic closure of $K$
with a fixed embedding $\bar{K} \rc \CC$
and all finite field extensions of K are regarded as subfields in $\bar{K}$\\
$\phi: H \to H$, the $\absfr{}$-linear frobenius morphism
for $H$ an $F$-(iso)crystal over $k$\\

Let $S = \Spec A \to T = \Spec B$ be a morphism of affine schemes
and let $X$ be an algebraic object over $S$.
We write $X \times_S T$ simply as $X_B$ or $X \tensor_A B$.

All deformations or lifts discussed in this paper
are supposed to be flat families.\\

\begin{center}
{\large {\bf I. Canonical and quasi-canonical lifts}}
\end{center}

\begin{center}
{\large 1. Basic properties}
\end{center}

\ni {\bf 1.1.}
Let $X$ be a K3 surface over $k$.
Assume that $X$ is ordinary,
i.e. the Newton polygon of the $F$-crystal $( H^2_{cris} (X/W), \phi)$
over $k$ coincides with its Hodge polygon.
In this case,
the $F$-crystal $H^2_{cris} (X/W)$ splits canonically as
\begin{equation*}
H^2_{cris} (X/W) = H^2(X, W{\cal O}) \oplus
				H^1(X, W\Omega^1) \oplus
				H^0(X, W\Omega^2),
\end{equation*}
and the absolute frobenius $\phi$ restricted to $H^{2-i} (X, W\Omega^i)$
is divisible by $p^i$ but not divisible by $p^{i+1}$.
The $W$-module $H^{2-i} (X, W\Omega^i)$ is free
of rank 1, 20, 1 for $i = 0, 1, 2$, respectively
(see \cite{Illusie}, p.653).\\

Let $\Phi$ be the formal Brauer group (\cite{AM}, p.109)
and let $\Psi$ be the enlarged formal Brauer group of $X$
(\cite{AM}, Proposition IV.1.8).
Then $\Phi$ is a $p$-divisible formal group and
its covariant Cartier module is canonically identified with
$H^2(X, W{\cal O})$ (\cite{AM}, Corollary II.4.3).
The group $\Phi$ is identified with the connected component
$\Psi^{\circ}$ of $\Psi$,
and $\Psi$ is a split extension
of an \'etale $p$-divisible group $\Psi^{\rm et}$ of rank $20$ by $\Phi$
\[ \Psi = \Phi \oplus \Psi^{\rm et}. \]\\

Let $\mathfrak{Art}_W$ be the category of artinian local $W$-algebras.
Let $\mathfrak{Def}_X$ be the deformation functor
that sends an object $S$ in $\mathfrak{Art}_W$
to the set of isomorphism classes of flat lifts of $X$ to $S$.
We extend $\mathfrak{Def}_X$ to formal lifts of $X$
over a complete discrete valuation $W$-algebra $(R, \mathfrak{m})$
by setting
\[ \mathfrak{Def}_X (R) := \projlim\ \mathfrak{Def}_X(R/\mathfrak{m}^n). \]
Similarly let $\mathfrak{Def}_{\Psi}$ be the deformation functor of $\Psi$.
For a lift $\mathfrak{X}$ of $X$ to $S$,
the associated enlarged formal Brauer group $\Psi_{\mathfrak{X}}$
of $\mathfrak{X}$ is a lift of $\Psi$ to $S$.
Thus we get a natural transformation
\[ \gamma: \mathfrak{Def}_X \to \mathfrak{Def}_{\Psi}. \]
It is known that $\gamma$ is an isomorphism of functors
when $X$ is ordinary
(\cite{Nyg-Tate}, Theorem 1.6).\\

\ni{\bf 1.2.}
Let $S$ be an object in $\mathfrak{Art}_W$
and let $\mathfrak{G}$ over $S$ be a lift of $\Psi$.
We have the connected-\'etale exact sequence
\[ 0 \to \mathfrak{G}^{\circ} \to \mathfrak{G} \to \mathfrak{G}^{\rm et} \to 0, \]
where $\mathfrak{G}^{\circ}$ is the connected part
and $\mathfrak{G}^{\rm et}$ is the maximal \'etale quotient
of $\mathfrak{G}$.
Since $\Phi$ and $\Psi^{\rm et}$ are rigid,
they have unique lifts to $S$, still call them $\Phi$ and $\Psi^{\rm et}$,
respectively.
Thus $\mathfrak{G}^{\circ} = \Phi$,
$\mathfrak{G}^{\rm et} = \Psi^{\rm et}$ over $S$
and $\mathfrak{G}$ corresponds to an element $[\mathfrak{G}]$
in ${\rm Ext} (\Psi^{\rm et}, \Phi)(S)$.
This gives an identification
\[ \mathfrak{Def}_{\Psi} = {\rm Ext} (\Psi^{\rm et}, \Phi), \]
and this provides a natural group structure on
$\mathfrak{Def}_X$.\\

\ni{\bf Definition 1.3} (\cite{LST}, \S 7).
Let $\Psi$ be the associated enlarged formal Brauer group
of an ordinary K3 surface $X$ over $k$.
Let $S$ be an artinian local $W$-algebra and
$\mathfrak{G}$ a lift of $\Psi$ to $S$.
The lift $\mathfrak{G}$ is called the {\em canonical lift} of $\Psi$
if the connected-\'etale sequence of $\mathfrak{G}$ splits over $S$, i.e.
if $\mathfrak{G}$ corresponds to the identity in the group
${\rm Ext} (\Psi^{\rm et}, \Phi)(S)$.
The lift $\mathfrak{G}$ is called a {\em quasi-canonical lift}
if the corresponding element in ${\rm Ext} (\Psi^{\rm et}, \Phi)(S)$
is torsion.
We also extend the definition to formal lifts
over a complete discrete valuation $W$-algebra.\\

\ni{\bf Definition 1.4.}
Let $X$ be an ordinary K3 surface over $k$
and $\Psi$ the enlarged formal Brauer group of $X$.
Let $S$ be an artinian local $W$-algebra.
Let $\mathfrak{X}$ be a lift of $X$ to $S$.
The lift $\mathfrak{X}$ is called the {\em canonical lift}
if the associated enlarged formal Brauer group $\Psi_\mathfrak{X}$
is the canonical lift of $\Psi$.
The lift $\mathfrak{X}$ is called a {\em quasi-canonical lift}
if $\Psi_\mathfrak{X}$ is a quasi-canonical lift of $\Psi$.
We also extend the definition to formal lifts
over a complete discrete valuation $W$-algebra.\\

\ni{\bf 1.5.}
Suppose $X$ is an ordinary K3 surface over $k$.
Now let $R$ be a complete discrete valuation $W$-algebra
with $[R:W] < \infty$.
Let $L$ be the field of fractions of $R$.
We fix a formal lift $\mathfrak{X}$ of $X$ to $R$.
If $\mathfrak{X}$ is algebraizable,
i.e. it is the formal completion of a (unique) scheme over $R$,
we then use the same $\mathfrak{X}$ to indicate the scheme
and set $X^{\circ} = \mathfrak{X} \tensor_R L$
to be the generic fiber of $\mathfrak{X}$.\\

\ni{\bf Lemma 1.6.}
{\em Keep the assumption as above.
If $\mathfrak{X}$ is a quasi-canonical lift,
then $\mathfrak{X}$ is algebraizable.}\\

\pf
This is a direct generalization of \cite{Nyg-Tate}, Proposition 1.8.
We keep the same notations there.
We may assume that $k$ is the residue field of $R$.
Let $\mathfrak{m}$ be the maximal ideal of $R$,
and $R_n = R/\mathfrak{m}^n$.
Write $\mathfrak{X} = \projlim X_n$,
where $X_n$ is a lift of $X = X_0$ to $R_n$
and the system $\{ X_n \}$ defines $\mathfrak{X}$.
By Grothendieck's existence theorem,
it suffices to find an ample line bundle on $X$
that can be lifted to all $X_n$.
Indeed, we will show that for any line bundle on $X$,
a certain power of it can be lifted to all $X_n$.

Let ${\bar X}_n = X_n \tensor_W W({\bar k})$.
Note $p$ is nilpotent on $X_n$ and $R_n$.
Consider the exact sequence of \'etale sheaves on $X_n$
\[ 1 \to 1 + \mathfrak{m} \mathcal{O}_{X_n} \to \mathcal{O}_{X_n}^*
	\to \mathcal{O}_X^* \to 1, \]
and the corresponding one on ${\bar X}_n$
\[ 1 \to 1 + \mathfrak{m} \mathcal{O}_{{\bar X}_n} \to
	\mathcal{O}_{{\bar X}_n}^* \to \mathcal{O}_{\bar X}^* \to 1. \]
We have a commutative diagram with injective vertical arrows
\[	\xymatrix{ H^1(X_n, \mathcal{O}_{X_n}^*) \ar[r] \ar@{_{(}->}[d] &
	H^1(X, \mathcal{O}_X^*) \ar[r] \ar@{_{(}->}[d] &
	H^2 (X_n, 1 + \mathfrak{m} \mathcal{O}_{X_n}) \ar@{_{(}->}[d] \\
	H^1({\bar X}_n, \mathcal{O}_{{\bar X}_n}^*) \ar[r] &
	H^1({\bar X}, \mathcal{O}_{\bar X}^*) \ar[r]^{\hspace{-20pt}\alpha} &
	H^2 ({\bar X}_n, 1 + \mathfrak{m} \mathcal{O}_{{\bar X}_n}). } \]
Thus it suffices to show that a certain multiple of the connecting map
$\alpha$ vanishes.
Henceforth we assume that $k$ is algebraically closed.

The map $\alpha$ factors through
\[ 	\xymatrix{ H^1(X, \mathcal{O}_X^*) \ar[r] &
	\projlim H^1(X, \mathcal{O}_X^*/ {\mathcal{O}_X^*}^{p^r})
	\ar[r]^{\projlim\ \beta_r} &
	H^2 (X_n, 1 + \mathfrak{m} \mathcal{O}_{X_n}) }, \]
where $\beta_r$ is the connected morphism
\[ H^1(X, \mathcal{O}_X^*/ {\mathcal{O}_X^*}^{p^r}) \to
	H^2 (X_n, 1 + \mathfrak{m} \mathcal{O}_{X_n}) \]
associated to the exact sequence of sheaves on $X_n$
\[ 1 \to 1 + \mathfrak{m} \mathcal{O}_{X_n} \to
	\mathcal{O}_{X_n}^* / {\mathcal{O}_{X_n}^*}^{p^r} \to
	\mathcal{O}_X^* / {\mathcal{O}_X^*}^{p^r} \to 1 \]
for $r \geq n$ (see \cite{Nyg-Tate}, p.218).
Let $\Psi_n$ be the enlarged formal Brauer group of $X_n$.
Let $E_{\Psi_n} \in {\rm Ext} (\Psi^{\rm et}, \Phi)(R_n)$
be the associated extension class of $\Psi_n$ and
\begin{eqnarray*}
	\delta: {\rm Ext} (\Psi^{\rm et}, \Phi)(R_n) &\isoto&
	\Hom_{\ZZ_p} \left( H^2_{fl} (X, \ZZ_p(1)), \Phi(R_n) \right) \\ &=&
	\Hom_{\ZZ_p}
	(\projlim H^1(X, \mathcal{O}_X^* / {\mathcal{O}_X^*}^{p^r} ),
	H^2(X, 1 + \mathfrak{m} \mathcal{O}_{X_n}) )
\end{eqnarray*}
be the isomorphism defined in
\cite{Nyg-Tate}, p.217 (cf. \cite{Messing}, p.180).
Then we have
$\projlim\ \beta_r (= \projlim ``p^r" ) = \delta( E_{\Psi_n} )$
(\cite{Nyg-Tate}, Proposition 1.7).
By assumption, the extension class $E_{\Psi_n}$ is torsion.
Hence a certain power of $\projlim\ \beta_r$ vanishes.
\qed\\

\ni{\bf Remark 1.7.}
If $k$ is of characteristic $p > 2$,
one can apply the {\em canonical coordinates}
on the deformation space of $X$ (\cite{Del-CC})
to provide another proof of Lemma 1.6.
Indeed, as in the proof above,
we can assume that $k$ is algebraically closed.
Then
there exists a system of canonical coordinates $(a, b_n, c, q_n)$
such that the set of formal functions $\{ q_n -1\}_{1 \leq n \leq 20}$
gives a system of coordinates on the universal deformation space
$S$ of an ordinary K3 surface $X$ over $k$
(\cite{Del-CC}, Th\'eor\`eme 2.1.7).
Thus
we have an isomorphism of formal schemes
\[ S \isom \Spf W[[q_1 - 1, \cdots, q_{20} - 1]]. \]
On the other hand,
we have $\Psi^{\circ} \isom \mu_{p^{\infty}}$ and
$\Psi^{\rm et} \isom (\QQ_p/\ZZ_p)^{\oplus 20}$ non-canonically
and
\[ S \isom {\rm Ext} (\Psi^{\rm et}, \Psi^{\circ}) \isom
	\widehat{\GG}_m^{\oplus 20} \]
is a formal torus (\cite{Messing}, Proposition A.2.5).
Under these isomorphisms,
$\{q_n\}$ is a basis of characters
\[ \chi: \widehat{\GG}_m^{\oplus 20}\ \to \ \widehat{\GG}_m \]
of the formal torus $\widehat{\GG}_m^{\oplus 20}$
(cf. \cite{Del-CC}, Th\'eor\`eme 2.1.14).

Suppose $\mathfrak{X}$ is a lift of $X$ over the ring of integers $R$
of a finite extension $L$ of $K$.
Then $\mathfrak{X}$ corresponds to an $R$-point $\lambda$ of $S$
\[ \begin{array}{cccc}
	\lambda: &W[[q_1 - 1, \dots, q_{20} - 1 ]] &\to& R \\
		& q_n &\mapsto& \lambda_n.
	\end{array} \]
If $\mathfrak{X}$ is a quasi-canonical lift,
the point $\lambda \in {\rm Ext} (\Psi^{\rm et}, \Psi^{\circ})(R)$ is torsion.
Thus $\lambda_n$ are roots of unit and are congruent to $1$
modulo the maximal ideal $\mathfrak{m}_R$ of $R$.
If $I$ is a non-trivial invertible sheaves on $X$,
then the deformation space $\mathfrak{Def}_{(X,I)}$
of the pair $(X, I)$
is a closed formal subscheme of $S$ defined by
\[ \prod_{1 \leq n \leq 20} q_n^{x_n} = 1\]
for $x_n \in \ZZ_p$ such that $c_1(I) = \sum x_n e_q^* b_n$,
where $e_q : W[[q_n]] \to W$ defined by $e_q(q_n) = 1$ for all $n$.
Note that the map $e_q$ gives the canonical lift of $X$
(\cite{Del-CC}, Th\'eor\`eme 2.2.2).
Since $c_1(I^{\tensor r}) = \sum r x_n e_q^* b_n$,
a power of $I$ can be lifted
to the quasi-canonical lift $\mathfrak{X}$.\\

\ni{\bf Lemma 1.8.} (Cf. \cite{Nyg-Tate}, Theorem 2.6.)
{\em Keep the assumption as in \S 1.5.
Then $\mathfrak{X}$ is a quasi-canincal lift of $X$
if and only if $\mathfrak{X}$ is projective and
the $\Gal(\bar{K}/L)$-module $H^2_{et}(X^{\circ}_{\bar K}, \QQ_p)$
attached to the generic fiber $X^{\circ} := \mathfrak{X}_L$
splits
\[ H^2_{et}(X^{\circ}_{\bar K}, \QQ_p) = V_0 \oplus V_1 \oplus V_2 \]
such that
\[ V_i \tensor_{\QQ_p} \CC_p \isom \CC_p(-i)^{\oplus h_i} \]
for $i = 0, 1, 2$ as $\Gal(\bar{K}/L)$-modules,
where $\CC_p$ is the completion of $\bar{K}$.
(In this case, $h_0 = h_2 = 1$ and $h_1 = 20$.)}
\\

\pf
Suppose $\mathfrak{X}$ is projective.
Let $\Psi$ be the enlarged formal Brauer group of $X$
and $\mathfrak{G}$ be the enlarged formal Brauer group of $\mathfrak{X}$.
Then the rational Tate module
\[ V(\mathfrak{G}) := \left( \projlim\ {\rm Ker} \left\{
	p^n: \mathfrak{G}({\bar K}) \to
	\mathfrak{G}({\bar K}) \right\} \right) \tensor_{\ZZ_p} \QQ_p \]
is a $\Gal({\bar K}/L)$-submodule of
$H^2_{et} (X^{\circ}_{\bar K}, \QQ_p(1))$
(\cite{AM}, Proposition IV.2.1 and Theorem IV.4.1).
Suppose that $H^2_{et} (X^{\circ}_{\bar K}, \QQ_p(1))$ splits
as in the Proposition.
Then the $p$-divisible group $\mathfrak{G}$
is isogenous to the direct product
$\Psi^{\circ}_R \times \Psi^{\rm et}_R$ over $R$.
Here,
$\Psi^{\circ}_R$ (resp. $\Psi^{\rm et}_R$) denotes
the unique lift of $\Psi^{\circ}$ (resp. $\Psi^{\rm et}$) to $R$.
Thus $\mathfrak{G}$ is a quasi-canonical lift of its special fiber $\Psi$
(\cite{LST}, \S7) and
therefore $\mathfrak{X}$ is a quasi-canonical lift of $X$.

Conversely,
if $\mathfrak{X}$ is quasi-canonical,
it is projective by Lemma 1.6
and the $p$-adic \'etale cohomology group
has the desired decomposition (see \cite{Nyg-Tate}, Theorem 2.6).
\qed\\

\ni{\bf Lemma 1.9.}
{\em Keep the assumption as in \S 1.5.
A lift $\mathfrak{X}$ over $R$ is quasi-canonical if and only if
the Hodge filtration ${\rm Fil}^{\bullet}$ on
$H_{dR}^2 (\mathfrak{X}/R) \tensor_R L$ coincides
with the one induced by the weight filtration
on $H^2_{cris}(X/W)$ via the natural isomorphism
\begin{equation*}
H^2_{dR} (\mathfrak{X}/R) \tensor_R L \isom
	H^2_{cris} (X/W) \tensor_W L
\tag{1.9.1}
\end{equation*}
(\cite{BO}, Corollary 2.5).
In particular,
\[ {\rm Fil}^2 \left( H^2_{dR} (\mathfrak{X}/R) \tensor_R L\right) 
	= H^0 (X, W\Omega^2) \tensor_W L. \]}
\\
\pf
If ${\rm Fil}^2 \left( H^2_{dR} (\mathfrak{X}/R) \tensor_R L\right) 
	= H^0 (X, W\Omega^2) \tensor_W L$,
then the whole Hodge filtration ${\rm Fil}^{\bullet}$ on
$H^2_{dR} (\mathfrak{X}/R) \tensor_R L$ coincides with
that induced by the weight filtration on $H^2_{cris}(X/W)$
since the identification (1.9.1) is compatible
with cup products.

By definition,
the lift $\mathfrak{X}$ is quasi-canonical if and only if
the enlarged formal Brauer group $\Psi_{\mathfrak{X}}$
attached to $\mathfrak{X}$ is quasi-canonical.
By applying the diagram in [NO], Theorem 3.20,
evaluating at $R$ and taking tensor with $\QQ$,
we see that $\Psi_{\mathfrak{X}}$ is quasi-canonical if and only if
${\rm Fil}^2 \left( H^2_{dR} (\mathfrak{X}/R) \tensor_R L\right) 
	= H^0 (X, W\Omega^2) \tensor_W L$.
\qed\\

\ni{\bf Remark 1.10.}
Assume that a formal lift $\mathfrak{X}$ of $X$ to $R$ is projective.
Then we have the $p$-adic Hodge comparison
\[ H^2_{et} (X^{\circ}_{\bar K}, \QQ_p) \tensor_{\QQ_p} {\bf B}_{cris}
	= H^2_{cris} (X/W) \tensor_W {\bf B}_{cris} \]
which extends to
\[ H^2_{et} (X^{\circ}_{\bar K}, \QQ_p) \tensor_{\QQ_p} {\bf B}_{dR}
	= H^2_{dR} (X^{\circ}/L) \tensor_L {\bf B}_{dR} \]
under (1.9.1) and a choice of a inclusion ${\bf B}_{cris} \rc {\bf B}_{dR}$
(\cite{Font-ss}, Th\'eor\`eme 6.1.4).
Then one sees that
a decomposition of the $\Gal(\bar{K}/L)$-representation
$H^2_{et} (X^{\circ}_{\bar K}, \QQ_p)$
corresponds to a decomposition of the filtered $\phi$-module
$H^2_{dR} (X^{\circ}/L)$.
Thus if the formal lift $\mathfrak{X}$ is projective,
the splitting of the $\Gal(\bar{K}/L)$-representation in Lemma 1.8 and
the coincidence of the filtrations in
Lemma 1.9 are equivalent to each other
by taking into account of the Tate twists.\\

\ni{\bf Definition 1.11.}
Keep the notations as in \S 1.5.
Assume that $k$ is a finite field and
the lift $\mathfrak{X}$ is projective.
We say that the frobenius morphism on $X$ is {\em liftable
to $H^2(\mathfrak{X}(\CC), \QQ)$}
if there exists an endomorphism $\pi$ on $H^2(\mathfrak{X}(\CC), \QQ)$
which satisfies the following two conditions: \\

i) The induced endomorphism $\pi \tensor_{\QQ} id$
on $H^2(\mathfrak{X}(\CC), \QQ) \tensor_{\QQ} \QQ_\ell$ 
for every prime $\ell \neq p$ is the geometric
frobenius on $H^2_{et} (X_{\bar k}, \QQ_\ell)$
via the natural identification
\begin{equation*}
H^2(\mathfrak{X}(\CC), \QQ_\ell) = H^2_{et}(\mathfrak{X}_{\CC}, \QQ_\ell)
	= H^2_{et}(\mathfrak{X}_{\bar K}, \QQ_\ell)\ \isoto\
	H^2_{et}(X_{\bar k}, \QQ_\ell)
\tag{1.11.1}
\end{equation*}
(hence $\pi \tensor_{\QQ} id$ commutes with Galois action).\\

ii) At the prime $p$,
we have $(\pi \tensor_{\QQ} id) \tensor_{\QQ_p} id
	= \phi^a \tensor_W id$
via the $p$-adic Hodge comparison
\begin{equation*}
\left( H^2(\mathfrak{X}(\CC), \QQ) \tensor_{\QQ} \QQ_p \right)
		\tensor_{\QQ_p} {\bf B}_{cris}
	= H^2_{et}(\mathfrak{X}_{\bar K}, \QQ_p)
		\tensor_{\QQ_p} {\bf B}_{cris}
	= H^2_{cris}(X/W) \tensor_W {\bf B}_{cris}
\tag{1.11.2}
\end{equation*}
(hence $(\pi \tensor_{\QQ} id) \tensor_{\QQ_p} id$
commutes with Galois action).

In this case, we say that the frobenius can be {\em lifted to
$H^2(\mathfrak{X}(\CC), \QQ)$},
and $\pi$ is called the {\em lifted frobenius}.\\

\ni{\bf Theorem 1.12.}
{\em Let $\mathfrak{X}$ be a formal lift of
an ordinary K3 surface $X$ over a finite field $k$ to
a complete discrete $W$-algebra $R$
with $[R:W] < \infty$.
Then $\mathfrak{X}$ is quasi-canonical if and only if
$\mathfrak{X}$ is projective and
a power of the frobenius morphism on $X$ is liftable to
$H^2 (\mathfrak{X}(\CC), \QQ)$.}\\

\pf
We may assume that $k$ is the residue field of $R$.
Suppose $\mathfrak{X}$ is a quasi-canonical lift.
Let $X^{\circ}$ be the generic fiber of $\mathfrak{X}$
over the field of fractions $L$ of $R$.
Then
after replacing $k$, $R$ and $L$ by compatible finite extensions,
there exist an abelian variety $A^{\circ}$ over $L$
with good reduction $\mathfrak{A}$ over $R$
and a $\ZZ$-algebra $C \rc \End_L A^{\circ}$
such that there exists an isomorphism of rational Hodge structures
\[ u: C^+P^2 (X^{\circ}(\CC), \QQ(1))\ \isoto\
	\End_C H^1 (A^{\circ}(\CC), \QQ) \]
which induces an isomorphism of $\Gal(\bar{K}/L)$-modules,
after extending the scalars to $\QQ_\ell$ for every prime $\ell$,
\[ u_l: C^+P^2_{et} (X^{\circ}_{\bar K}, \QQ_\ell(1))\ \isoto\
	\End_C H^1_{et} (A^{\circ}_{\bar K}, \QQ_\ell), \]
where $P^2_{\bullet}(\bullet)$ is the primitive part
in the corresponding cohomology group
with respect to a chosen polarization and
$C^+$ is the even Clifford algebra
(\cite{Del-K3}, Proposition 6.5 and \S 6.6).
Let $A$ over $k$ be the special fiber of $\mathfrak{A}$.
Then $A$ is an ordinary abelian variety and
$\mathfrak{A}$ is isogenous to the canonical lift of $A$
(\cite{Nyg-Tate}, Proposition 2.5 and Corollary 2.8).
Let $\pi_A \in (\End_R \mathfrak{A}) \tensor_{\ZZ} \QQ$
that lifts the frobenius morphism of $A$.
Via $u$, the map $\pi_A$ induces an endomorphism $\pi$
on $H^2(X^{\circ}(\CC), \QQ)$ which satisfies (1.11.1)
(cf. \cite{Nyg-Tate}, Lemma 3.2).

At the prime $p$,
by applying $u_p$, tensoring with ${\bf B}_{cris}$ over $\QQ_p$
and taking the $\Gal(\bar{K}/L)$-invariants,
we get an isomrphism of $F$-isocrystals
\[ u_{cris}: C^+ \left( P^2_{cris} (X/W) \tensor_W K(1) \right)\
	\isoto\ \End_C \left( H^1_{cris}(A/W) \tensor_W K \right) \]
(\cite{Font-ss}, Th\'eor\`eme 6.1.4).
Since $\pi_A$ is the lifted frobenius morphism on $A$,
the induced map of $\pi_A$ on $H^1_{cris}(A/W) \tensor K$
is just $\phi^a$,
where $\phi$ is the absolute frobenius on the crystalline cohomology
and $p^a$ is the number of elements in $k$.
Then it is easy to see that
the induced map of $\pi$ on $H^2_{cris} (X/W)$
is the geometric frobenius morphism
since $u_{cris}$ transfers the absolute frobenius
to the absolute frobenius.

Conversely,
by a base change of $k$,
we assume that the frobenius endomorphism of $X$ over $k$
can be lifted to $\pi$ on $H^2(X^{\circ}(\CC), \QQ)$.
Suppose the field $k$ has $p^a$ elements.
Since $X$ is ordinary,
the characteristic polynomial $f_{\phi}(x)$
of the absolute frobenius $\phi$ on $H^2_{cris}(X/W)$ decomposes
$f_{\phi}(x) = f_{\phi}^{(0)}(x) f_{\phi}^{(1)}(x) f_{\phi}^{(2)}(x)$
according to the valuation of its roots.
(Although $f_{\phi}(x)$ is not well-defined,
the valuations of its roots, the $a$-th powers of the roots,
and the decomposition of the crystal
$H^2_{cris}(X/W)$ are well-defined.)
Now the roots of the characteristic polynomial $f_{\pi}(x)$
on $H^2(X^{\circ}(\CC), \QQ)$
are $a$-th powers of the roots of $f_{\phi}(x)$.
Thus $f_{\pi}(x)$ decomposes as
$f_{\pi}(x) = f_{\pi}^{(0)}(x) f_{\pi}^{(1)}(x) f_{\pi}^{(2)}(x)$ mod $p$.
By Gauss Lemma,
$f_{\pi}(x)$ has a decomposition over $\ZZ_p$
that lifts the decomposition mod $p$.
Thus $H^2_{et} (X^{\circ}_{\bar K}, \QQ_p)$
has a decomposition as in Lemma 1.8 and
consequently $\mathfrak{X}$ is a quasi-canonical lift.
\qed\\

\begin{center}
{\large 2. Hodge structure}
\end{center}

\ni{\bf 2.1.}
Back to the situation as in \S 1.5.
Then Theorem 1.12 says that
there exists an endomorphism $\pi^n$
on $H^2 (X^{\circ}(\CC), \QQ)$
which lifts the $n$-th iteration of the frobenius endomorphism
of $X$ relative to $k$ for some positive integer $n$.\\

\ni{\bf Lemma 2.2.}
{\em Keep the assumption as above.
With respect to the embedding $L \rc \bar{K} \rc \CC$,
the Hodge structure $H^2(X^{\circ}(\CC), \QQ)$
of the complex projective K3 surface $X^{\circ} \tensor_L \CC$
is determined by the lifted frobenius $\pi^n$
and the map $\pi^n$ is a Hodge cycle.}\\

\pf
We need to figure out the 1-dimensional sub-$\CC$-vector space
\[ H^0 (X^{\circ}_{\CC}, \Omega^2_{X_{\CC}}) =
	{\rm Fil}^1 \left( H^2_{dR} (X^{\circ}_{\CC}/\CC)(1) \right) \rc
	H^2_{dR} (X^{\circ}_{\CC}/\CC)(1) = H^2(X^{\circ}(\CC), \CC)(1). \]
We consider $\bar{\QQ}$ as a subfield of $\CC$.
Thus the embedding $\bar{K} \rc \CC$ provides
an embedding $\bar{\QQ} \rc \bar{K}$
and $\bar{\QQ}$ is equipped with the induced $p$-adic valuation.

Let $H = H^2(X^{\circ}(\CC), \QQ(1))$ and let
\[ \begin{array}{ccc}
	H_{>0} &=&
		\left( H^0(X, W\Omega^2) \tensor_W K(1) \right)
			\tensor_K \CC \\
	H_0 &=&
		\left( H^1(X, W\Omega^1) \tensor_W K(1) \right)
			\tensor_K \CC \\
	H_{<0} &=&
		\left( H^2(X, W\mathcal{O}) \tensor_W K(1) \right)
			\tensor_K \CC. \\
\end{array} \]
Then on $H_{\CC} = H^2_{dR}(X^{\circ}/L)(1) \tensor_L \CC$,
we have 
\[ H_{\CC} = H_{>0} \oplus H_0 \oplus H_{<0}. \]
Since $\mathfrak{X}$ is a quasi-canonical lift,
the filtration is given by
${\rm Fil}^1 H_{\CC} = H_{>0}$,
which is 1-dimensional over $\CC$,
and $\pi$ acts on it through the first factor
$H^0(X, W\Omega^2) \tensor_W K(1)$ of $H_{>0}$
as multiplication by a constant $q\epsilon$,
where $\epsilon$ is a $p$-adic unit.
Thus ${\rm Fil}^1 H_{\CC}$ is the unique complex 1-dimensional
eigenspace of $\pi$ in $H_{\CC}$
attached to the unique eigenvalue with positive $p$-adic valuation.
\qed\\

\ni{\bf Lemma 2.3.}
{\em Keep the assumption as in \S 2.1.
Then the frobenius morphism on $X$ relative to $k$
can be lifted to an endomorphism $\pi$ on $H^2(X^{\circ}(\CC), \QQ)$
and the map $\pi$ is a Hodge cycle.}\\

\pf
We will abuse notations to
let $\pi$ be the geometric frobenius
on any cohomology of $X$.
Take an integer $n$ such that
the iterated frobenius $\pi^n$ can be lifted to
$H^2(X^{\circ}(\CC), \QQ)$, still call it $\pi^n$,
and such that $\pi^n = p^a$ on the algebraic cycle classes in
$H^2(X^{\circ}(\CC), \QQ)$ for a suitable integer $a$.
By Lemma 2.2, we can
regard $\pi^n$ as a class $[\pi^n]$ in
$\End H^2(X^{\circ}(\CC), \QQ) \rc
	H^4(X^{\circ}(\CC) \times X^{\circ}(\CC), \QQ(2))$,
where $\End H^2(X^{\circ}(\CC), \QQ)$
is the algebra of endomorphisms of the rational Hodge structure
$H^2(X^{\circ}(\CC), \QQ)$.

Pick a prime $\ell \neq p$.
Then the class $[\pi^n]$ in
$H^4(X^{\circ}(\CC) \times X^{\circ}(\CC), \QQ_\ell(2)) =
H^4_{et}(X_{\bar k} \times X_{\bar k}, \QQ_\ell(2))$
is an algebraic cycle class of $X_{\bar k}$
since it is represented by the graph $\Gamma_{\pi^n}$
of the endomorphism $\pi^n$ on $X$.
We know that
the set of algebraic cycle classes
(modulo homological equivalence with $\QQ$ coefficients)
$C^2(X_{\bar k} \times X_{\bar k}, \QQ) \rc
H^4_{et}(X_{\bar k} \times X_{\bar k}, \QQ_\ell(2))$
are $\QQ$-spanned by
$pr_i^*[X]$, ${\rm NS}(X_{\bar k}) \times {\rm NS}(X_{\bar k})$
and the graphs of the iterations of $\pi^n$,
where $pr_i$ are the projections $X \times X \to X$
to the $i$-th factor for $i = 1,2$
(\cite{Zarhin-cycle}, Corollary 6.1.1).
Notice in this case,
${\rm NS}(X^{\circ}_{\CC})_{\QQ} \to {\rm NS}(X_{\bar k})_{\QQ}$
is an isomorphism (\cite{Nyg-Tate}, Theorem 3.3)
and we can regard
${\rm NS}(X^{\circ}_{\CC}) \times {\rm NS}(X^{\circ}_{\CC})$
as a subspace in $\End H^2(X^{\circ}(\CC), \QQ)$.
Thus the class $[\pi]$ in $\End H^2(X^{\circ}(\CC), \QQ)$
can be obtained by a $\QQ$-combination of
$[\pi^n], [\pi^{2n}] = [\pi^n] \circ [\pi^n], \cdots$ and
elements in
${\rm NS}(X^{\circ}_{\CC}) \times {\rm NS}(X^{\circ}_{\CC})$.
\qed\\

\ni{\bf 2.4.}
Now let $Z$ be a K3 surface over the complex numbers $\CC$.
Let $\SS = {\rm Res}_{\CC/\RR} \GG_m$ be the restriction of scalars
from $\CC$ to the real numbers $\RR$
of the multiplicative group $\GG_m$.
Let ${\rm NS}(Z) \rc H^2(Z, \ZZ(1))$ be the N\'eron-Severi group of $Z$.
The cup product restricted to ${\rm NS}(Z)$ is non-degenerate.
Let ${\rm NS}(Z)^{\perp}$ be the orthogonal complement of
${\rm NS}(Z)$ in $H^2(Z, \ZZ(1))$ with respect to the cup product pairing.
Then ${\rm NS}(Z)^{\perp}$ is a sub-integral Hodge structure
of $H^2(Z, \ZZ(1))$.
Let ${\sf M}(Z) = {\rm NS}(Z)^{\perp} \tensor_{\ZZ} \QQ$.
Then there is an orthogonal decomposition
of the rational Hodge structure
\[ H^2 (Z, \QQ(1)) = {\rm NS}(Z)_{\QQ} \oplus {\sf M}(Z). \]
\\

Let
\[ h: \SS \to {\rm GL} \left( {\sf M}(Z)_{\RR} \right) \]
be the associated homomorphism that induces the Hodge structure
on ${\sf M}(Z)$.
Let ${\sf Hdg}_Z$ be the Hodge group of $h$, i.e.
it is the smallest algebraic group over $\QQ$
in ${\rm GL} ({\sf M}(Z))$ such that
its $\RR$-valued points contain the image of
$\{ z \in \CC^* | z\bar{z} = 1 \} \rc \SS(\RR) = \CC^*$ under $h$.
(In this case,
it equals the Mumford-Tate group of ${\sf M}(Z)$
since the weight of the Hodge structure ${\sf M}(Z)$ is zero.)
Let ${\sf E} = \End_{{\sf Hdg}_Z} {\sf M}(Z)$.
Thus ${\sf M}(Z)$ is a ${\sf E}$-vector space
with an ${\sf E}$-linear ${\sf Hdg}_Z$-action.\\

With these notations,
one knows that
the ${\sf Hdg}_Z$-module ${\sf M}(Z)$
is irreducible (\cite{Zarhin-Hodge}, Theorem 1.4.1).
The $\QQ$-algebra ${\sf E}$ is either a totally real number field
or a CM number field (\cite{Zarhin-Hodge}, Theorem 1.5.1).
Finally
the group ${\sf Hdg}_Z$ is commutative if and only if
${\sf E}$ is a CM field and
$\dim_{\sf E} {\sf M}(Z) = 1$.
In this case, ${\sf Hdg}_Z \rc {\rm Res}_{{\sf E}/\QQ} \GG_m$
and $\dim_{\QQ} {\sf Hdg}_Z = \frac{1}{2} \dim_{\QQ} {\sf E}$
(\cite{Zarhin-Hodge}, Theorem 2.2.1, and 2.3.1).\\

We now consider the case
where the complex K3 surface $Z$
is from a quasi-canonical lift $\mathfrak{X}$
of an ordinary K3 surface $X$ as in \S 2.1.\\

\ni{\bf Theorem 2.5.}
{\em Assume $k$ is a finite field of $q$ elements
and $X$ over $k$ is an ordinary K3 surface.
Let $L$ be a finite extension of $K$ with ring of integers $R$.
Let $\mathfrak{X}$ be a quasi-canonical lift of $X$
and $X^{\circ} = \mathfrak{X} \tensor_R L$.
As above,
Let ${\sf Hdg}_X$ be the Hodge group
of the transcendental part ${\sf M}(X^{\circ}_{\CC})$
of $H^2(X^{\circ}(\CC), \QQ(1))$
and ${\sf E} = \End_{{\sf Hdg}_X} {\sf M}(X^{\circ}_{\CC})$.
Let ${\sf G} = \overline{<\pi / q>}$ be the Zariski closure over $\QQ$
of the cyclic group $< \pi / q >$ generated by $\pi/q$
in ${\rm GL} \left( {\sf M}(X^{\circ}_{\CC}) \right)$.
Then we have
${\sf E} = \QQ[\pi]$ and ${\sf Hdg}_X = {\sf G}$.}\\

\pf
Since $\pi$ respects the Hodge structure on ${\sf M}(X^{\circ}_{\CC})$,
we have $\pi \in {\sf E}$.
Since the linear map $\pi$ is irreducible on ${\sf M}(X^{\circ}_{\CC})$
(\cite{Zarhin-cycle}, Theorem 1.1),
we have $\QQ[\pi] = {\sf E}$
as ${\sf M}(X^{\circ}_{\CC})$ is an ${\sf E}$-vector space.

Let $2n = \dim_{\QQ} {\sf E}$.
Given an embedding ${\sf E} \rc \CC$ sending $\pi$ to $\alpha \in \CC$,
it induces a character of the torus $\sf{G}$
by sending the generator $\pi/q$ to $\alpha/q$.
In this way, we get $n$ independent characters
(\cite{Zarhin-cycle}, Theorem 1.2).
Thus $\dim_{\QQ} \sf{G} \geq n$.
On the other hand,
every element in $\sf{G}(\QQ)$ has norm one
with respect to any embedding ${\sf E} \rc \CC$.
Thus $\dim_{\QQ} {\sf G} \leq n$.
Therefore $\dim_{\QQ} {\sf G} = n$
and it is the maximal compact subtorus of
${\rm Res}_{{\sf E}/\QQ} \GG_m$.
Consequently, $\sf{Hdg}_X = \sf{G}$ by the description
of the shape of $\sf{Hdg}_X$ and dimension counting
(\cite{Zarhin-Hodge}, Remark 2.3.2).
\qed\\

\ni{\bf Remark 2.6.}
i)
As in \cite{Zarhin-Hodge}, \S 1.1,
we extend the action of ${\sf Hdg}_X$ to the full Hodge structure
$H^2 (X^{\circ}(\CC), \QQ(1))$ by trivial action on
${\rm NS}(X^{\circ}_{\CC})_{\QQ}$.
Then the action of $\pi/q$,
regarded as an element in ${\sf Hdg}_X(\QQ)$,
is not the same as its action regarded as the lifted frobenius
since the lifted frobenius may act non-trivially on
${\rm NS}(X^{\circ}_{\CC})_{\QQ}$.\\

ii)
Since ${\sf Hdg}_X$ is a reductive Lie group,
it is the stabilizer group of a finite set $S$ of Hodge cycles
in the tensor algebra of ${\sf M}(X^{\circ}_{\CC})$.
We know that
every Hodge cycle in the tensors generated by ${\sf M}(X^{\circ}_{\CC})$
is absolutely Hodge (\cite{LNM900}, Proposition II.6.25)
hence we may assume that $(\pi/q)^n$ fixes every elements in $S$
for some positive integer $n$.
Thus $(\pi/q)^n \in {\sf Hdg}_X(\QQ)$
since ${\sf Hdg}_X$ is the largest $\QQ$-subgroup of
${\rm GL}({\sf M}(X^{\circ}_{\CC}))$ that fixes all Hodge cycles.
This provides a direct way to build an inclusion
${\sf G} \rc {\sf Hdg}_X$.\\

\begin{center}
{\large {\bf II. Applications}}
\end{center}

\begin{center}
{\large 3. Torelli theorem}
\end{center}

\ni{\bf 3.1.}
Let $X$ be an ordinary K3 surface over a finite field $k$.
Let $\mathfrak{X}$ over $W$ be the canonical lift of $X$.
Write $X^{\circ} = \mathfrak{X}_K$.
Let $T(X) = H^2(X^{\circ}(\CC), \ZZ)$ and
$V(X) = T(X) \tensor_{\ZZ} \QQ = H^2 (X^{\circ}(\CC), \QQ)$.
The $\ZZ$-module $T(X)$ is equipped with a quadratic form
induced from the cup product pairing.
It is isomorphic to
\[ \UU \oplus \UU \oplus \UU \oplus (-\EE_8) \oplus (-\EE_8), \]
where $\UU$ is the hyperbolic plane and
$\EE_8$ is the unique positive definite even unimodule of rank 8.
For any rational prime $\ell$,
let $T_\ell(X) = T(X) \tensor_{\ZZ} \ZZ_\ell$
and $V_\ell(X) = V(X) \tensor_{\QQ} \QQ_\ell$.
Then Lemma 2.3 says that
the frobenius morphism on $X$ over $k$
can be lifted as an endomorphism $\pi$ on $V(X)$.\\

\ni{\bf Proposition 3.2.}
{\em Assume $p > 3$.
The lifted frobenius $\pi$ on $V(X)$ preserves the $\ZZ$-lattice $T(X)$.}\\

\pf
We show that for any prime $\ell$,
the lifted frobenius $\pi$ preserves the $\ZZ_{\ell}$-lattice
$T_{\ell}(X)$ in $V_{\ell}(X)$.
For $\ell \neq p$,
the map $\pi \tensor_{\QQ} id$
on $H^2(X^{\circ}(\CC), \QQ_{\ell})$
is the geometric frobenius on $H^2_{et} (X_{\bar k}, \QQ_{\ell})$.
Thus it preserves the lattice
$H^2_{et}(X_{\bar k}, \ZZ_{\ell}) = T_{\ell}(X)$.

At the prime $p$,
since $\mathfrak{X}$ over $W$ is the canonical lift of $X$,
by Lemma 1.9,
the Hodge filtration on $H^2_{dR}(\mathfrak{X}/W)$
coincides with the weight filtration on $H^2_{cris}(X/W)$.
Thus the geometric frobenius is an endomorphism
on the filtered $\phi$-module $H^2_{cris}(X/W)$.
Therefore $\pi$ preserves the $\ZZ_p$-lattice
$H^2_{et}(X^{\circ}_{\bar K}, \ZZ_p) = T_p(X)$
in $H^2_{et}(X^{\circ}_{\bar K}, \QQ_p)$ if $p > 3$ (\cite{FM}, Remark 6.4).
\qed\\

\ni{\bf Theorem 3.3.}
(Cf. \cite{Nyg-Torelli}, Theorem 2.1.)
{\em Assume $p > 3$.
Let $(X, I)$ and $(Y, J)$ be polarized ordinary K3 surfaces
over a finite field $k$
and $(T(X), \pi_X)$ and $(T(Y), \pi_Y)$ be the associated
quadratic $\ZZ$-modules with lifted frobenius endomorphisms
of $X$ and $Y$, respectively.
Let $\xi \in T(X)$ and $\eta \in T(Y)$ be the (twisted)
first chern classes of $I$ and $J$, respectively.
Suppose there exists an isometry
\[ f^*: T(Y) \to T(X) \]
which commutes with the frobenius
and sends $\eta$ to $\xi$ .
Then there exists an isomorphism $f: (X, I) \to (Y, J)$ over $k$
such that $f^*_\ell := f^* \tensor id: T_\ell(Y) \to T_\ell(X)$
is induced from $f$ on the \'etale cohomology under
$T_\ell(X) = H^2_{et} (X^{\circ}_{\CC}, \ZZ_\ell)
	= H^2_{et}(X_{\bar k}, \ZZ_\ell)$
for all primes $\ell \neq p$.}\\

\pf
Let $(\mathfrak{X, I})$ and $(\mathfrak{Y, J})$
be the canonical lifts of $(X, I)$ and $(Y, J)$ respectively
(cf. Remark 1.7;
\cite{Nyg-Tate}, Proposition 1.8).
Since the map $f^*: T(Y) \to T(X)$
commutes with the frobenius $\pi_Y$ and $\pi_X$,
it respects Hodge structures by Lemma 2.3.
Thus by the classical Torelli theorem for K3 surfaces over $\CC$
(\cite{BPV}, Theorem VIII.11.1),
there exists an isomorphism
$f_{\CC}: (\mathfrak{X, I})_{\CC} \to (\mathfrak{Y, J})_{\CC}$
that induces $f^*$.
Choose an isomorphism $\overline{K} \isom \CC$,
we may assume that the map $f_{\CC}$ is defined over $\overline{K}$
\[ f_{\overline{K}}: (\mathfrak{X, I})_{\overline{K}}\ \to\
	(\mathfrak{Y, J})_{\overline{K}}. \]
Again,
since the map $f_{\overline{K}}$ commutes with frobenius,
it is defined over a valuation field $L \rc {\overline{K}}$
with ring of integers $R$ such that the residue field of $R$ is $k$.
Thus we get an isomorphism
\[ f_L: (\mathfrak{X, I})_L\ \to\ (\mathfrak{Y, J})_L \]
such that $f_{\overline{K}} = f_L \tensor_L \overline{K}$.

Since the morphism $f_L$ is an isomorphism
between the generic fibers of
polarized smooth non-ruled schemes
$(\mathfrak{X, I})$ and $(\mathfrak{Y, J})$
over the ring of integers $R$ of $L$,
it
can be extended uniquely to an isomorphism
\[ f_R: (\mathfrak{X, I})_R\ \to\
	(\mathfrak{Y, J})_R \]
(\cite{MM}, Theorem 2).
In particular,
we get the desired isomorphism
\[ f: X \to Y \]
over $k$ by taking reduction of $f_R$.
\qed\\

\ni{\bf Remark 3.4.}
It is easy to see that at the prime $\ell = p$,
the map $f^*_{cris}$
defined by the commutativity of the diagram
\[	\xymatrix{
	H^2_{et} (\mathfrak{Y}_{\bar K}, \QQ_p) \tensor_{\QQ_p} {\bf B}_{cris}
		\ar[r]^{\sim} \ar[d]_{f^*} &
	H^2_{cris} (Y/W) \tensor_W {\bf B}_{cris}
	\ar[d]^{f^*_{cris}} \\
	H^2_{et} (\mathfrak{X}_{\bar K}, \QQ_p) \tensor_{\QQ_p} {\bf B}_{cris}
	\ar[r]^{\sim} &
	H^2_{cris} (X/W) \tensor_W {\bf B}_{cris}
	}\]
is induced by the morphism $f: X \to Y$ on the crystalline cohomology.\\

\begin{center}
{\large 4. Hodge conjeture}
\end{center}

\ni{\bf 4.1.}
Suppose we are in the situation of \S 1.5.
As in \S 2.4,
we let ${\sf M}(X^{\circ}_{\CC})$ be the transcendental part
of $H^2(X^{\circ}(\CC), \QQ(1))$.
The subspace ${\sf M}(X^{\circ}_{\CC})$
is self-dual via the symmetric product induced by the cup product
\[ H^2(X^{\circ}(\CC), \QQ(1)) \times H^2(X^{\circ}(\CC), \QQ(1))\
	\to\ H^4(X^{\circ}(\CC), \QQ(2)) = \QQ. \]
Under this self-duality,
we regard the lifted frobenius $\pi$ on ${\sf M}(X^{\circ}_{\CC})$
as an element in ${\sf M}(X^{\circ}_{\CC})^{\tensor 2}$.\\

\ni{\bf Proposition 4.2.}
{\em Keep the assumption as above.
Then the subspace of Hodge cycles in the tensor algebra of
${\sf M}(X^{\circ}_{\CC})$
is generated by the graphs $\Gamma_n$
of iterations $\pi^n$
of the lifted frobenius $\pi$ in ${\sf M}(X^{\circ}_{\CC})^{\tensor 2}$.}\\

\pf
Let $C^n \rc {\sf M}(X^{\circ}_{\CC})^{\tensor 2n}$
be the subspace of Hodge cycles.
Then by Theorem 2.5,
the space $C^n \tensor_{\QQ} \QQ_\ell \rc
	{\sf M}(X^{\circ}_{\CC})^{\tensor 2n} \tensor_{\QQ} \QQ_{\ell}$
coincides with Tate cycles in
${\sf M}(X^{\circ}_{\CC})^{\tensor 2n} \tensor_{\QQ} \QQ_{\ell}$
regarded as a subspace in
$H^2(X^{\circ}(\CC), \QQ_{\ell}(1))^{\tensor 2n} =
	H^2_{et}(X_{\bar k}, \QQ_{\ell}(1))^{\tensor 2n}$.
Since the $\QQ_{\ell}$-space of Tate cycles
in the tensor algebra of
${\sf M}(X^{\circ}_{\CC}) \tensor_{\QQ} \QQ_{\ell}$
regarded as a subspace in
$H^2_{et}(X_{\bar k}, \QQ_{\ell}(1))$
is generated by $\Gamma_n$
(\cite{Zarhin-Tate}, Corollary 6.1.1),
the same holds for the space of Hodge cycles
in the tensor algebra of ${\sf M}(X^{\circ}_{\CC})$.
\qed\\

\ni{\bf A.} Kummer surfaces.\\

\ni{\bf 4.3.}
Let $A$ be an ordinary abelian variety
over a finite field $k$
with $q = p^a$ elements.
Assume that the frobenius morphism acts trivially
on the 2-torsion points $A[2]$ of $A$.
Let $L$ be a finite extension of K in ${\bar K}$
with ring of integers $R$.
Let $\mathfrak{A}$ be a lift of $A$ over $R$
and $A^{\circ}$ be the generic fiber of $\mathfrak{A}$.
Let $\delta: \widetilde{\mathfrak{A}} \to \mathfrak{A}$
be the blow-up of $\mathfrak{A}$
along the 2-torsion points $\mathfrak{A}[2]$.
Then the special fiber $\widetilde{A}$ of $\widetilde{\mathfrak{A}}$
is the blow-up of $A$ along $A[2]$
and the generic fiber $\widetilde{A}^{\circ}$ of $\widetilde{\mathfrak{A}}$
is the blow-up of $A^{\circ}$ along $A^{\circ}[2]$.
We have
\[ \delta^* : H^2(A^{\circ}(\CC), \ZZ) \to
	H^2(\widetilde{A}^{\circ}(\CC), \ZZ) \isom 
	H^2(A^{\circ}(\CC), \ZZ) \oplus \ZZ(-1)^{\oplus 16}, \]
where the subgroup $\ZZ(-1)^{\oplus 16}$ is generated by the sixteen
(-1)-curves from the blow-up.
Extending the scalars,
$\delta^*$ induces a homomorphism $\delta_\ell^*$
of $\Gal({\bar K}/L)$-modules
\[ \delta_\ell^*: H_{et}^2(A_{\bar K}, \ZZ_\ell) \to
	H_{et}^2(\widetilde{A}_{\bar K}, \ZZ_\ell) \isom 
	H_{et}^2(A_{\bar K}, \ZZ_\ell) \oplus \ZZ_\ell(-1)^{\oplus 16}, \]
for every prime number $\ell$. \\

Let $\iota: \mathfrak{A} \to \mathfrak{A}$ be the involution.
Then $\iota$ extends to an involution
$\iota: \widetilde{\mathfrak{A}} \to \widetilde{\mathfrak{A}}$.
Let $\mathfrak{X}$ be the Kummer surface associated to $\mathfrak{A}$
and $\epsilon: \widetilde{\mathfrak{A}} \to \mathfrak{X}$
the quotient map (which is generically two-to-one).
Let $X$ be the special fiber and $X^{\circ}$ the generic fiber
of $\mathfrak{X}$.
Then $X$ (resp. $X^{\circ}$) is the Kummer surface
associated to $A$ (resp. $A^{\circ}$).
We have the following diagram:
\[ \xymatrix{
	& \widetilde{\mathfrak{A}} \ar[dl]_{\delta} \ar[dr]^{\epsilon} & \\
	\mathfrak{A} \ar@{-->}[rr]^{2:1} & & \mathfrak{X}. }
		\tag{4.3.1} \]
In this case
\[ \epsilon^*: H^2(X^{\circ}, \QQ) \to
	H^2(\widetilde{A}^{\circ}, \QQ) \]
is an isomorphism of rational Hodge structures
and $\epsilon^*$ is compatible with cup-product pairings.
Extending the scalars, $\epsilon^*$
induces an isomorphism $\epsilon^*_\ell$
of $\Gal({\bar K}/L)$-modules
\[ \epsilon^*_\ell: H^2_{et} (X_{\bar K}, \QQ_\ell) \to
	H^2_{et} (\widetilde{A}_{\bar K}, \QQ_\ell) \]
for every prime number $\ell$.\\

\ni{\bf Lemma 4.4.}
{\em Keep assumptions as above. Then

{\rm a)} The Kummer K3 surface $X$ is an ordinary K3 surface.

{\rm b)} The formal scheme $\mathfrak{X}$ is a quasi-canonical lift
of $X$ if $\mathfrak{A}$ is a quasi-canonical lift of $A$.

{\rm c)} If $\mathfrak{X}$ is a quasi-canonical lift of $X$,
then $\mathfrak{A}$ is a quasi-canonical lift of $A$.

{\rm d)} If $L = K$ and $\mathfrak{A}$ is the canonical lift of $A$,
then $\mathfrak{X}$ is the canonical lift of $X$.}\\

\pf
a)
Fix an embedding $\bar{\QQ} \rc \bar{K}$.
After a base change,
we may assume that the whole Kummer construction (4.3.1) is defined
over $R$ with residue field $k$.
Take a prime $\ell \neq p$.
We have an isomorphism of $\Gal({\bar k}/k)$-modules
\[ \epsilon_{\ell}^*: H^2_{et} (X_{\bar k}, \QQ_\ell)\ \isoto\
	H^2_{et} (\widetilde{A}_{\bar k}, \QQ_\ell). \]
Thus the characteristic polynomial $f(X; x)$
of the geometric frobenius on $H^2_{et}(X_{\bar k}, \QQ_\ell)$
equals that $f(\widetilde{A}; x)$ on
$H^2_{et} (\widetilde{A}_{\bar k}, \QQ_\ell)$.
Since
\[ H^2_{et}(\widetilde{A}_{\bar k}, \QQ_\ell) =
	H^2_{et}(A_{\bar k}, \QQ_\ell) \oplus \QQ_\ell(-1)^{\oplus 16} \]
and $A$ is ordinary,
there exists a unique root of $f(\widetilde{A}; x)$
which is a $p$-adic unit.
Therefore $X$ is ordinary.\\							

b)
If $\mathfrak{A}$ is a quasi-canonical lift of $A$,
then
\[ H^1_{et} (A^{\circ}_{\bar K}, \QQ_p) = H_0 \oplus H_1 \]
with
\[ H_i \tensor_{\QQ_p} \CC_p \isom \CC_p(-i)^{\oplus 2} \]
as $\Gal({\bar K}/L)$-modules.
Thus
\[ H^2_{et} (A^{\circ}_{\bar K}, \QQ_p)
	= \left( H_0 \wedge H_0 \right)
	\oplus \left( H_0 \tensor H_1 \right)
	\oplus \left( H_1 \wedge H_1 \right) \]
and therefore
\[ H^2_{et} (X^{\circ}_{\bar K}, \QQ_p) = V_0 \oplus V_1 \oplus V_2 \]
where
\[ \begin{array}{ccl}
	V_0 &=& \left(\epsilon_p^*\right)^{-1} \left( H_0 \wedge H_0 \right) \\
	V_1 &=& \left(\epsilon_p^*\right)^{-1} \left( H_0 \tensor H_1
		\oplus \QQ_p(-1)^{\oplus 16} \right) \\
	V_2 &=& \left(\epsilon_p^*\right)^{-1} \left( H_1 \wedge H_1 \right).
\end{array} \]
Since
\[ V_i \tensor \CC_p \isom \CC_p(-i)^{\oplus h_i}, \]
$\mathfrak{X}$ is a quasi-canonical lift by Lemma 1.8.\\

c) Assume that $\mathfrak{X}$ is a quasi-canonical lift.
Let $\mathfrak{B}$ be the Kuga-Satake abelian scheme
over the ring of integers $R'$
of a finite extension $L'$ of $L$ insider ${\bar K}$.
Then after replacing $L'$ by another finite extension,
$\mathfrak{A} \tensor_R R'$ is isogenous to a factor of $\mathfrak{B}$
(\cite{Sk}, Theorem 2).
Since $\mathfrak{B}$ is a quasi-canonical lift of its special fiber
(cf. \cite{Nyg-Tate}, Corollary 2.8),
\[ H^1_{et}(\mathfrak{B}_{\bar K}, \QQ_p) = H_0 \oplus H_1 \]
with
\[ H_i \tensor_{\QQ_p} \CC_p \isom \CC_p(-i)^{\oplus h}. \]
Consequently $H^1_{et}(A^{\circ}_{\bar K}, \QQ_p)$
has a similar decomposition as a $\Gal({\bar K}/L')$-module.
Thus $\mathfrak{A} \tensor_R R'$ is a quasi-canonical lift of $A$.
Since the deformation of $A$ is pro-representable,
it is left exact.
Therefore $\mathfrak{A}$ is a quasi-canonical lift of $A$ to $R$.\\

d)
As in c), since the deformation of $X$ is pro-representable,
we may replace $k$ by its algebraic closure.
Thus we assume that $k$ is algebraically closed.
Since $\mathfrak{X}$ is defined over $W$
and is a quasi-canonical lift by part b),
it is the canonical lift of $X$
as there is only one torsion element
in the group $1 + pW$ of principal units. (cf. Remark 1.7).
\qed

\ni{\bf Theorem 4.5.}
{\em Let $A$ be a complex abelian surface of CM type.
Let $X$ be the associated Kummer surface.
Then the Hodge conjecture is true for any self product
$X \times \cdots \times X$ of $X$.}\\

\pf
If $A$ is isogenous to a product of two elliptic curves,
Then the statement is true
since the Hodge conjecture is known for products of
arbitrary elliptic curves.
Henceforth, assume that $A$ is simple.
Since $A$ is of CM type,
it is define over a number field $F$.
Write $E := (\End_F A) \tensor_{\ZZ} \QQ$,
which is a CM field with $[E:\QQ] = 4$.

For a non-archimedean place $v$ of $F$,
let $F_v$ be the completion of $F$ at $v$.
Let $A_v = A \tensor_F F_v$,
$\mathfrak{A}_v$ its N\'eron model and
$\bar{A}_v$ the special fiber of $\mathfrak{A}_v$.
Let $\pi_v$ be the frobenius endomorphism of $\bar{A}_v$
relative to the residue field $k_v$ at $v$.
Enlarging $F$ if necessary,
there exists a non-archimedean place $v$
such that the reduction $\bar{A}_v$ is a simple ordinary abelian variety
over $k_v$
and the characteristic polynomial of $\pi_v$ has no multiple root
(\cite{Tankeev}, Lemma 7.5).
Thus we have
$(\End_{k_v} \bar{A}_v) \tensor_{\ZZ} \QQ = \QQ(\pi_v)$
and $[\QQ(\pi_v) : \QQ] = 4$ (\cite{Tate}, Theorem 2).
Consequently, $E = \QQ(\pi_v)$
since $\End_F A \ \to\ \End_{k_v} \bar{A}_v$
is injective.
Therefore $\pi_v$ can be lifted to an isogeny on $A_v$,
which implies that $\mathfrak{A}_v$
is a quasi-canonical lift of $\bar{A}_v$.
Thus by Lemma 4.4,
$X$ can be regarded as the generic fiber
of a quasi-canonical lift of the Kummer surface
associated to $\bar{A}_v$
and the statement then follows by Proposition 4.2.
\qed\\

\ni{\bf B.} Weighted K3 surfaces.\\

\ni{\bf 4.6.}
Let $X$ be one of the weighted K3 surfaces considered in
\cite{Yo} and \cite{Yui}.
It is of one of the following two types:\\

\ni a).
The K3 surface $X$ is the resolution of singularities
of a certain hypersurface of Fermat type
\[ x_0^{n_0} + x_1^{n_1} + x_2^{n_2} + x_3^{n_3} \]
in the weighted projective 3-space $\mathbb{P}(m_0, m_1, m_2, m_3)$
(see \cite{Yui}, \S 2).
In this case,
we consider $X$ as a variety defined over a number field $F$
containing $E := \QQ(\zeta_n)$
where $n$ is the least common multiple of $n_i$, $i = 0, \dots, 3$
and $\zeta_n$ is a primitive $n$-th root of unity.
Then $E$ acts on any cohomology group of $X$
and this action commutes with the extra structure
attached to that cohomology theory (Hodge structure, Galois action, \dots).
There are fourteen K3 surfaces of this type.\\

\ni b).
The K3 surface $X$ is the resolution of singularities
of a certain quasi-diagonal hypersurface 
\[ x_0^{n_0} + x_0 x_1^{n_1} + x_2^{n_2} + x_3^{n_3} \]
in the weighted projective 3-space $\mathbb{P}(m_0, m_1, m_2, m_3)$
(see \cite{Yui}, \S 3).
In this case,
we consider $X$ as a variety defined over a number field $F$
containing $E := \QQ(\zeta_n)$
where $n$ is the least common multiple of $n_1, n_2$, and $n_3$.
Then $E$ acts on any cohomology group of $X$
and this action commutes with the extra structure
attached to that cohomology theory.
There are eighty five K3 surfaces of this type.\\

In both cases,
we can find an open set $S$ of the ring of integers of $F$
such that $X$ can be extended to $\mathfrak{X}$ over $S$.
For a finite place $v$ in $S$,
let $F_v$ be the completion of $F$ at $v$,
$R_v$ be the ring of integers of $F_v$ and
$k_v$ the residue field of $R_v$.
Set $X_v = X \tensor_F F_v$,
$\mathfrak{X}_v = \mathfrak{X} \tensor R_v$
and $\bar{X}_v$ the reduction of $\mathfrak{X}_v$.\\

\ni{\bf Theorem 4.7.}
{\em Let $X$ be as above.
Then there exists a finite place $v$ of $F$ such that
$\bar{X}_v$ is an ordinary K3 surface over the finite field $k_v$,
and $\mathfrak{X}_v$ is a quasi-canonical lift of $\bar{X}_v$.
Thus the Hodge conjecture is true for any self product
$X \times \cdots \times X$ of $X$.}\\

\pf
In fact,
there exists a positive integer $m$
such that if $v$ is a place over a rational prime $p$
with $p \equiv 1$ (mod $m$),
then the reduction $\bar{X}_v$ is ordinary
(\cite{GY}, Proposition 3.8).
Moreover
$\rk {\rm NS}(X_{\CC}) = 22 - [E: \QQ]
	= \rk {\rm NS}(\bar{X}_v \tensor_{k_v} \bar{k}_v)$
(\cite{LNM900}, Proposition I.7.6 and
\cite{Yui}, Proposition 6.2).
Similarly to \S 2.4,
we pick a prime $\ell \neq p$ and let
${\sf M}_{\ell}(\bar{X}_v)$
be the orthogonal complement of
${\rm NS}(\bar{X}_v \tensor_{k_v} \bar{k}_v)$ in
$H^2_{et} (\bar{X}_v \tensor_{k_v} \bar{k}_v, \QQ_{\ell}(1))$.
Then both $\zeta_n$ and
the frobenius endomorphism $\pi$ on $\bar{X}_v$
can be regarded as algebraic cycles in
${\sf M}_{\ell}(\bar{X}_v)^{\tensor 2}$.
Since $\dim_{\QQ} \QQ(\zeta_n) =
	\dim_{\QQ_{\ell}} {\sf M}_{\ell}(\bar{X}_v) =
	\dim_{\QQ} \QQ(\pi)$
(\cite{Zarhin-cycle}, Theorem 1.1),
we must have $\QQ(\zeta_n) = \QQ(\pi)$
(\cite{Zarhin-Tate}, Corollary 6.1.1).
Thus $\pi$ can be lifted to $H^2(X(\CC), \QQ)$
and $X_v$ is a quasi-canonical lift of $\bar{X}_v$
by Theorem 1.12.
\qed\\

\ni{\bf Remark 4.8.}
After this paper has been written,
the author learned that in the preprint \cite{Mari},
it has been shown that
the Hodge conjecture holds for arbitrary self products
of any complex projective K3 surface whose
Hodge endomorphism algebra on its transcendental lattice
is a CM field.
In particular this covers our examples here.
The proof in \cite{Mari} uses a result of
Mukai on the study of moduli of projective bundles over K3 surfaces
and isogenies between K3 surfaces,
which is different from our reduction approach.\\


\end{document}